\documentclass[12pt]{article}
\usepackage{amssymb,amsmath,amsfonts,amsthm,amscd,latexsym,indentfirst,geometry}
\def\id{\textrm{id}}

\def\char{\textrm{char}\,}
\def\Var{\mathrm{Var}}
\def\Span{\mathrm{span}}
\def\Aut{\mathrm{Aut}}
\def\Imm{\mathrm{Im}\,}

\def\qed{\hfill $\square$\vspace{10pt}}
\def\qedd{\hfill $\square$}

\begin{document}

\begin{center}
\noindent {\Large \bf Rota---Baxter Operators on Quadratic Algebras}
\vspace{0.5cm}

\noindent {\large Pilar Benito, Vsevolod Gubarev and Alexander Pozhidaev}
\\[3mm]

\end{center}

\parbox[c]{14cm}{{\bf Abstract.}  \it 
We prove that all Rota---Baxter operators on a quadratic division
algebra are trivial. For nonzero weight, we state that all
Rota---Baxter operators on the simple odd-dimensional Jordan algebra
of bilinear form are projections on a subalgebra along another one.
For weight zero, we find a connection between the Rota---Baxter
operators and the solutions to the alternative Yang---Baxter equation
on the Cayley---Dickson algebra. We also investigate the Rota---Baxter
operators on the matrix algebras of order two, the Grassmann algebra
of plane, and the Kaplansky superalgebra.} \\*[3mm]

\parbox[c]{14cm}{{\bf Mathematics Subject Classification.}  16T25, 17A45, 17C50.}\\*[2mm]

 \parbox[c]{14cm}{{\bf Keywords.}  \it
Rota---Baxter operator, Yang---Baxter equation,
quadratic algebra, matrix algebra, Grassmann algebra, Jordan algebra
of bilinear form, Kaplansky superalgebra.}

\section{Introduction}

Given an algebra $A$ and a scalar $\Delta$ in a field $F$, a linear
operator $R\colon A\rightarrow A$ is called a {\it Rota---Baxter
operator} ({\it  RB-operator}, shortly) on $A$ of weight $\Delta$ if
the following identity
\begin{equation}\label{RB}
R(x)R(y) = R( R(x)y + xR(y) + \Delta xy )
\end{equation}
holds for all $x,y\in A$. The algebra $A$ is called the {\it
Rota---Baxter algebra} ({\it RB-algebra}).

The Rota---Baxter algebras were introduced by Glen Baxter in 1960
\cite{GB}, and then they were popularized  by G.-C.~Rota and his
school \cite{RG1,RG2}. The linear operators with the property
\eqref{RB} were independently introduced in the context of Lie
algebras by A.~A.~Belavin and V.~G.~Drinfeld in 1982
\cite{BelaDrin82} and by M.~A.~Semenov-Tyan-Shansky in 1983
\cite{Semenov83}. These operators were connected with the so-called
$R$-matrices, which are solutions to the classical Yang---Baxter
equation. Recently, some applications of the Rota---Baxter algebras
were found in such areas as the quantum field theory, the
Yang---Baxter equations, the cross products, the operads, the Hopf
algebras, the combinatorics and the number theory (some references
may be found, for example, in \cite{GEF}).

In 2000, M.~Aguiar established  a connection between the
Rota---Baxter algebras and the dendriform algebras. He showed that a
Rota---Baxter algebra of weight  $\Delta = 0$ possesses the structure
of a dendriform algebra. Later on, a connection with the dendriform
trialgebras was established \cite{EF}. Some functors between the
categories of the Rota---Baxter algebras and the dendriform
dialgebras (trialgebras) were investigated in \cite{GEF}.

In the present article we are interested in the study
(classification) of the structures of Rota---Baxter algebras on some
well-known simple (super)algebras. The investigations of this type
previously were carried out for the direct sum of the complex
numbers field in~\cite{Braga}, and the simple three-dimensional Lie
algebra $\mathrm{sl}_2(\mathbb{C})$ \cite{sl2,sl2-0}. In
\cite{Goncharov1}, M.~E.~Goncharov considered the structures of
bialgebra on an arbitrary simple finite-dimensional algebra $A$ over
a field of characteristic zero with a semisimple Drinfeld double. He
proved that these structures induce on $A$ Rota---Baxter operators of
nonzero weight. Also, for simple Lie algebras and some non
skew-symmetric solutions to the classical Yang---Baxter equations, he
constructed Rota---Baxter operators of nonzero weight. As a
corollary, he constructed Rota---Baxter operators of nonzero weight
on the simple non-Lie Malcev algebra.

Some of the results of the present article were proved by the
authors independently. Preliminary Section 2 consists of the results
of Pilar Benito (PB) and Vsevolod Gubarev (VG). The results of
Section 3 were obtained by PB and Alexander Pozhidaev (AP). The
results of Subsection 4.1 were proved by VG and of 4.3  --  by PB
and VG. The results of Subsections 4.4, 4.5, 5.2, and 5.3 were
obtained by VG, and they are actually some applications of the
technique developed in Section~3. Theorem 6 (5.1), which is a
reproof of \cite{Mat2}, was proved by VG. The results of Subsection
5.4 were obtained by~AP.

In what follows, the characteristic of the main field $F$ is different from two.

\section{Preliminaries}

By the {\it trivial} RB-operators of weight $\Delta$
 we mean the zero operator and $-\Delta\id$,
where $\id$ denotes the identity map.

Consider some well-known examples of RB-operators (see, e.g.,
\cite{GuoMonograph}).

{\bf Example 1}.
Given an algebra $A$ of continuous functions on $\mathbb{R}$,
an integration operator
$R(f)(x) = \int\limits_{0}^x f(t)\,dt$
is an RB-operator on $A$ of weight zero.

{\bf Example 2}.
Given an invertible derivation $d$ of an algebra $A$,
$d^{-1}$ is an RB-operator on $A$ of weight zero.

{\bf Example 3}.
Let $A = \{(a_1,a_2,\ldots,a_k,\ldots)\mid a_i\in \Bbbk\}$ be a countable sum of a field~$\Bbbk$
with the termwise addition, multiplication and scalar product.
An operator $R$ defined as
$R(a_1,a_2,\ldots,a_k,\ldots)
 = \big(a_1,a_1+a_2,\ldots,\sum\limits_{i=1}^k a_i,\ldots\big)$
is an RB-operator on $A$ of weight $-1$.

Note that the algebra $A$ from Example 3 is not simple as an algebra but it is simple as an RB-algebra.
Also this example may be generalized for an arbitrary variety of algebras.

{\bf Statement 1} \cite{GuoMonograph}. \it\
Let  $P$ be an RB-operator of weight $\Delta$. Then

a$)$ the operator $-P-\Delta\id$ is an RB-operator of weight $\Delta$,

b$)$ the operator $\Delta^{-1}P$ is an RB-operator of weight 1,
provided that $\Delta\neq0$.\rm\qed

Let $A$ be an algebra. In what follows we fix the notation $\phi$
for the map defined on the set of all RB-operators on $A$ as
 $\phi(P) = -P-\Delta(P)\id$. It is clear that $\phi^2$ coincides
 with the identity map.

{\bf Statement 2}. \it
Let $P$ be an RB-operator  of weight $\Delta$ on an algebra $A$,
and let $\psi\in\Aut(A)$. Then $P^{(\psi)} = \psi^{-1}P\psi$
is an RB-operator on $A$ of weight $\Delta$. \rm

{\it Proof} is straightforward.\qed

{\bf Statement 3} \cite{GuoMonograph}. \it
Assume that an algebra $A$ is splitted as a vector space
into the direct sum of two subalgebras $A_1$ and $A_2$.
An operator $P$ defined by the rule
\begin{equation}\label{Split}
P(x_1+x_2) = -\Delta x_2,\quad x_1\in A_1,x_2\in A_2,
\end{equation}
is an RB-operator on $A$ of weight $\Delta$.\rm\qed

The RB-operator from Statement 3 is a {\it splitting} RB-operator with
respect to the subalgebras $A_1$ and $A_2$. In \cite{Jian}, such
RB-operator is called a quasi-idempotent operator.

{\bf Remark 1}.
Let $P$ be a splitting RB-operator on an algebra $A$ of weight $\Delta$ with respect to subalgebras $A_1,A_2$.
Then $\phi(P)$ is an RB-operator of weight $\Delta$,
$$
\phi(P)(x_1+x_2) = -\Delta x_1, \quad x_1\in A_1,x_2\in A_2,
$$
and $\phi(P)$ is a splitting RB-operator with respect to the same
subalgebras $A_1$ and $A_2$.

{\bf Remark 2}.
The set of all splitting RB-operators on an algebra $A$ is in
bijective correspondence with all decompositions of $A$
 into the direct sum of two subalgebras.

{\bf Example 4}. \cite{Jian}
Let $A$ be an associative algebra, and let $e\in A$ be an element such that $e^2 = - \lambda e$, $\lambda\in F$.
A linear map $l_e\colon x\to ex$ is an RB-operator of weight $\lambda$ satisfying $R^2 + \lambda R = 0$.
If $\lambda\neq0$ then $l_e$ is a splitting RB-operator on $A$
with respect to the subalgebras $A_1 = (1-e)A$ and $A_2 = eA$,
and the decomposition $A = A_1\oplus A_2$ is exactly a Pierce one.

In an alternative algebra $A$ with an element $e$ such that
$e^2 = - \lambda e$, $\lambda\in F$, the operator $l_e$ is an RB-operator
if $e$ lies in the associative or commutative center of $A$.
It follows easily using the identities of alternative algebras \cite{Nearly}.

{\bf Example 5}. In \cite{Poj}, there were described all possible
linear Rota---Baxter structures on a 0-dialgebra with a bar-unit.

{\bf Example 6}.
In \cite{BBGN2011}, it was proved that every RB-algebra of weight $\Delta$
in the variety $\Var$ with respect to the operations
$$
x\succ y = R(x)y,\quad
x\prec y = xR(y),\quad
x\cdot y = \Delta xy
$$
is a post-$\Var$-algebra.

In \cite{Embedding}, given a post-$\Var$-algebra $A$, its enveloping RB-algebra $B$
of weight $\Delta$ in the variety $\Var$ was constructed.
By the construction, $B = A\oplus A'$,
where $A'$ is a copy of $A$ as a vector space,
and the RB-operator $R$ was defined as follows:
$R(a')=\Delta a$, $R(a)=-\Delta a$, $a\in A$.
From the definition we have $A_1=\ker R=\Span\{a+a'\mid a\in A\}$,
$A_2 = R(B) = A$, and $R$ is a splitting RB-operator on $B$ with respect to
$A_1$ and $A_2$.
So, given a post-$\Var$-algebra $A$, there exists an enveloping
algebra $B$ with a splitting RB-operator $R$ of weight~1.

{\bf Lemma 1} \cite{Gub2017}.\it\
Let $A$ be a unital algebra, and let $P$ be an RB-operator on $A$ of weight~$\Delta$.

a$)$ If $P(P(x)+\Delta x) = 0$ then $P$ is splitting.

b$)$ If $\Delta\neq0$ and $P(1)\in F$ then $P$ is splitting.

c$)$ If $\Delta = 0$ then $1\not\in \Imm P$. Moreover, if $A$ is a
simple finite-dimensional algebra, $\dim A>1$, then $\dim \ker P\geq2$.

d$)$ If $\Delta = 0$ and $P(1)\in F$ then $P(1) = 0$,
$P^2 = 0$, and $\Imm P\subset \ker P$. \rm

{\it  Proof}.
a) Show that $A = \ker P\oplus P(A)$ as the direct sum of vector spaces.
On the contrary, assume that there exists a nonzero $x\in\ker P\cap P(A)$.
Then, $x = P(y)$ and $P(x) = P^2(y) = 0$. By the hopothesis,
$x = P(y) = -(1/\Delta) P^2(y) = 0$, a contradiction.

By \eqref{RB}, $\ker P$ and $P(A)$ are some subalgebras of $A$.
From $P(P(x)+\Delta x) = 0$, we have that the restriction of $P$ on
$P(A)$ is equal to $-\Delta\id$, and $P(\ker P) = 0$.

b) By \eqref{RB} for $x=y=1$ we have $P(1)\in\{0,-\Delta\}$.
It suffices to consider only the case $P(1)=0$. Indeed, if $P(1)=-\Delta$,
by Statement 1, we can study an RB-operator $\phi(P)$ of the same weight,
and $\phi(P)(1)=0$. By Remark 1, we are done.

By \eqref{RB}, for $x\in A$, we have \begin{equation}\label{Lemma}
0 = P(1)P(x) = P(P(1)x+1\cdot P(x)+\Delta x) = P(P(x)+\Delta x).
\end{equation}
So, we apply a).

c) Suppose $R(x) = 1$ for some $x\in A$. By \eqref{RB},
$1 = R(x)R(x) = 2R(x) = 2$, a contradiction.

Let $A$ be a simple finite-dimensional algebra, $\dim A = n$.
By a), $\dim \Imm P\leq n-1$. Assume that $\dim \Imm P = n-1$.
By \eqref{RB}, $\ker P$ is an $\Imm P$-bimodule.
Since $A = \Span\{1,\Imm P\}$, $\ker P$ is  a proper ideal of $A$,
a contradiction with the simplicity of~$A$.

d) By c), $P(1) = 0$. Other assertions follow from
$$\qquad\qquad\qquad\qquad
0 = P(1)P(x) = P(P(1)x + 1\cdot P(x)) = P(P(x)).
\qquad\qquad\qquad\qquad\square
$$

{\bf Lemma 2}.\it\
Let $A$ be an algebra, and let $R$ be an RB-operator on $A$ of weight zero.

a$)$ A nonzero element $e\in A$ such that $e^2 = \alpha e$, $\alpha\in F^*$,
could not be an eigenvector of $R$ with nonzero eigenvalue.

b$)$ If $A$ is a unital finite-dimensional algebra, $\Imm (R)$ is abelian,
and $F$ is algebraically closed, then $R$ is nilpotent.\rm

{\it  Proof}.
a) If $R(e) = ke$ with $k\in F^*$ then
$$
\alpha k^2 e = k^2 e^2
 = R(e)R(e)
 = R(R(e)e + eR(e))
 = 2k R(e^2)
 = 2\alpha k R(e) = 2\alpha k^2 e,
$$
a contradiction.

b) Suppose that $v$ is an eigenvector of $R$ with nonzero eigenvalue $k$. We have
$$
0 = R(1)R(v)
  = R(R(1)v + R(v)) = R((1/k) R(1)R(v)+R(v))
  = R^2(v) = k^2 v,
$$
a contradiction.\qedd

\section{Quadratic algebras}

Let $A$ be a quadratic algebra, i.e., every element $x\in A$
satisfies the equation
\begin{equation}\label{Quadratic}
x^2-t(x)x+n(x)1 = 0,
\end{equation}
where 1 is a unit of $A$, the {\it trace} $t(x)$ is linear on $A$,
and the {\it norm} $n(x)$ is quadratic on~$A$ \cite{Nearly}.

Putting $f(x,y)=n(x+y)-n(x)-n(y)$, we get
\begin{equation}\label{CircInQuad}
x\circ y = t(x)y+t(y)x-f(x,y)1.
\end{equation}
We have $A = F1\oplus A_0$, where $A_0 = \{x\in A\mid t(x) = 0\}$.

Let $R$ be an RB-operator on $A$ of weight $\Delta$.
Setting $x = y$ in \eqref{RB} and applying \eqref{Quadratic}, we infer that
\begin{equation}\label{NormInQuad}
- n(R(x))1 = R(t(x)R(x)-f(x,R(x))+\Delta t(x)x-\Delta n(x)).
\end{equation}
Taking $x\in A_0$ in \eqref{NormInQuad} we obtain
\begin{equation}\label{FinalQuad}
n(R(x))1 = (f\bigl(x,R(x)\bigr)+\Delta n(x))R(1).
\end{equation}

Applying Lemma 1, we arrive at the following statement.

{\bf Lemma 3}.\it\
Let $A$ be a quadratic algebra with an RB-operator $R$ of weight $\Delta$.

a$)$ If $R(1) = 0$ or $R(1)\not\in F$,
then $n(R(x)) = 0$ and $R(x)(R(x) - t(R(x))1) = 0$ for all $x\in A_0$.

b$)$ For $\Delta = 0$, either $R(1) = 0$ or
$n(x+R(x)) = n(x)$ for all $x\in A_0$.

c$)$ For $\Delta \neq 0$, if $n(R(x))\neq 0$ for some $x\in A_0$,
then $R$ is splitting.\rm\qed

{\bf Theorem 1}.\it\
All RB-operators on a quadratic division algebra are trivial.\rm

{\it  Proof}.
If a quadratic division algebra $A$ coincides with $F$
then the statement is obvious. Let $\dim_F(A)\geq2$.

If $R(1) \in F$ then we have $R^2 = -\Delta R$ by Lemma~1.
For $\Delta = 0$, by Lemma~1\,d), $R(1) = 0$
and by Lemma~3\,a), $R(x)(R(x) - t(R(x))1) = 0$ for all $x\in A$.
Since $A$ has no zero divisors, $R(x)\in F$ for all $x\in A$.
By \eqref{RB}, $R = 0$.

For $\Delta\neq0$, by Lemma~1\,b),
$R$ is splitting with respect to some subalgebras $A_1$ and $A_2$,
i.e., $R(A_1) = 0$, and $R$ is equal to  $-\Delta\id$ on $A_2$.
Up to $\phi$ we have $1\in A_1$. For each $x \in A_2$, we have $x\in R(A_2)$;
by \eqref{Quadratic} and Lemma~3\,a),
we obtain $x(x - t(x)1) = 0$. As $x\not\in F$, $x = 0$. So $R = 0$.

Let $R(1)\not\in F$. By Lemma~3\,a),
$R(x)(R(x) - t(R(x))1) = 0$ for all $x\in A_0$.
When $\Delta = 0$, $R(x)\in F$ for all $x\in A_0$,
and  $R(x) = 0$ for every $x\in A_0$ by \eqref{RB}.
So, $\Imm R$ is the linear span of $R(1)$.
By \eqref{RB}, $R(1)^2 = \alpha R(1)$ for some $\alpha\in F$.
We have either $R(1) = 0$ and $R = 0$ or $R(1)\in F$, a contradiction.
For $\Delta \neq 0$,
we have $R(x)\in F$ for all $x\in A_0$
and by~\eqref{RB} we infer that $A_0$ is a proper ideal of $A$,
a contradiction with the divisibility.\qed

{\bf Corollary 1}.\it\
Given a quadratic division algebra $A$,
there are no representations of $A$ as a sum
$($as vector spaces$)$ of its proper subalgebras.\rm

{\it  Proof}.
Assume that $A$ is equal to $A_1\oplus A_2$,
where $A_1,A_2$ are some subalgebras of~$A$.
Hence, by Statement 3, there exist nontrivial RB-operators on $A$
of nonzero weight. By Theorem~1, we arrive at a contradiction.\qed

{\bf Lemma 4}.\it\
Let $A$ be a quadratic commutative algebra.
Then the RB-operators $R$ of weight $0$ on $A$ such that
$R(1) = 0$ are in one-to-one correspondence with the linear maps $R$ on
$A$ such that $R(1) = 0$, $\Imm R\subseteq \ker R\cap \ker n$.\rm

{\it  Proof}.
Let $R$ be an RB-operator of weight 0 on $A$ such that $R(1) = 0$.
By Lemma~1\,d), $R^2 = 0$. By \eqref{FinalQuad}, $n(R(A)) = 0$. Thus,
$\Imm R\subseteq \ker n$, and $\ker n\subseteq  \ker R$.

Conversely, let $R$ be a map on $A$ as above. Then $R^2 = 0$, and
$n(R(A)) = 0$. By \eqref{Quadratic}, $R(x)R(x) = t(R(x))R(x)$.
By \eqref{CircInQuad},
$$
R(x\circ R(x)) = R(t(x)R(x) + t(R(x))x) = t(R(x))R(x) = R(x)R(x).
$$
Thus, $A$ is a Rota---Baxter algebra of weight zero. \qed

Let $A$ be an algebra over a field $F$, let $S$ be a subalgebra of~$A$,
let $I$ be a subspace of $A$ such that  $SI+IS\subseteq I$, and
let $D$ be a nondegenerate derivation from $S$ to $A$ modulo $I$,
(i.e., $D(xy)-D(x)y-xD(y)\in I$ for all $x,y\in S$) with the property
$A = D(S)\oplus I$. In this case we say that $(S,I,D)$ is
an $RB$-{\it triple} on $A$.
Denote the space of all derivations from $S$ to $A$ modulo $I$
by $\mathrm{Der}_F(S,I,A)$.

{\bf Lemma 5}.\it\  Let $A$ be an algebra over a field $F$. Then the
RB-operators of weight~0 on $A$ are in one-to-one
correspondence with the $RB$-triples on $A$. \rm

{\it  Proof}.
Let $R$ be an RB-operator of weight~0 on $A$.
Put $I = \ker R$ and $S = \Imm R$. Choose a basis for $I$
and complete it to a basis of $A$ by some $a_j\in A$, $j\in J$
for some set of indexes $J$ such that $S = \Span\{R(a_j)\mid j\in J\}$.
Put $A_0 = \Span\{a_j\mid j\in J\}$. Then $A = I\oplus A_0$.
Define a linear map $D\colon S\to A$ by the rule $D(R(a_j))=a_j$ for all $j\in J$.
Then $A = D(S)\oplus I$. Note that if $a\in A$ then $a = i+a_0$ for
some uniquely defined $i\in I$, $a_0\in A_0$; therefore,
$D(R(x))\equiv x\,({\rm mod}\,I)$. Take arbitrary $a_i,a_j\in A_0$
and put $s_1=R(a_i)$, $s_2=R(a_j)$. Then by \eqref{RB} with $\Delta=0$ we have
$s_1s_2=R(s_1a_j+a_is_2)$, and
$$
D(s_1s_2)=D(R(s_1D(s_2)+D(s_1)s_2)
 \equiv s_1D(s_2)+D(s_1)s_2\, ({\rm mod}\, I).
$$
It is easy to see that $(S,I,D)$ is an $RB$-triple on $A$.

Conversely, let $(S,I,D)$ be an $RB$-triple on $A$.
Define an operator $R$ on $A$ by the rule
$$
\ker R=I,\quad R(D(s))=s,\ s\in S.
$$
If either $x\in I$ or $y\in I$ then \eqref{RB} holds. Take
$x=D(s_1)$, $y=D(s_2)$ for arbitrary $s_1,s_2\in S$. Then
\begin{gather*}
R(x)R(y)=s_1s_2,\\
R(R(x)y+xR(y))
 = R(s_1D(s_2)+D(s_1)s_2)
 = R(D(s_1s_2)) = s_1s_2,
\end{gather*}
and \eqref{RB} holds again.\qed

{\bf Corollary 2}.\it\
Let $\mathcal{V}$ be a variety of algebras over a field $F$.
Let $A$ be a $\mathcal{V}$-algebra, and let $V$ be an
$A$-module in the sense of Eilenberg.
Assume that there exists a nondegenerate derivation $D$ from $A$
into $B = A\oplus V$ modulo $V$ such that $B = D(A)\oplus V$.
Then $(A,V,D)$ is an RB-triple on $B$.\rm

{\it  Proof}.
By the definition of module in the sense of Eilenberg, we
have $A\leq B$, $AV+VA\subseteq V$. Now, apply Lemma~5.\qed

{\bf Remark 3}.
The hypotheses of Corollary~2 hold if $D$ is a nondegenerate
derivation~$D$ of $A$ such that $D(A) = A$.

{\bf Example 7}.
Consider the Lie algebra $\mathrm{sl}_2(\mathbb{C})$ with the standard basis $h,e,f$.
Put $S = \Span\{h\}$, $D = \mathrm{ad}\,(e+f)$, $I = \Span\{h,e\}$.
Then the operator $R$ such that $R(f) = h/2$ and $R(I) = 0$ gives
the RB-operator on $\mathrm{sl}_2(\mathbb{C})$ of weight zero.
It is exactly the case (R5) \cite{Kolesnikov} from six possible variants of RB-operators
on $\mathrm{sl}_2(\mathbb{C})$ of weight zero.

{\bf Example 8}.
Let $A$ be an algebra. Assume that $S$ is an abelian subalgebra of $A$,
$A = S\oplus I$, and $S$ acts on $I$.
(For example, one may consider a Lie algebra and its Cartan
subalgebra as $S$.) Then every nondegenerate mapping on $S$ with the kernel $I$
determines an RB-operator on $A$ of weight zero.

{\bf Example 9}.
Consider a semisimple finite-dimensional Lie algebra $L$ over a
 field $F$ of characte\-ris\-tic~0.
Assume that there are some nonzero roots
$\alpha,\beta$ such that $\beta+\alpha$
 belongs to the set $\Gamma$ of nonzero roots of $L$ but
 $\beta-\alpha\not\in \Gamma$. Take $h$ in the Cartan subalgebra $H$
of $L$ such that $\beta(h)\neq 0$. Put
$S = \Span\{h,e_{\alpha}\}$,
$I = H\oplus\sum\limits_{\gamma\in\Gamma\setminus\{\beta,\alpha+\beta\}}
\Span \{e_\gamma\}$, and $D = \mathrm{ad} (e_{\beta})$.
Consider the operator $R$ on $L$ such that $R(e_{\beta})=-\beta(h)^{-1}h$,
$R(e_{\alpha+\beta})=c_{\alpha,\beta}^{-1}e_{\alpha}$, where
$[e_{\alpha,\beta}]=c_{\alpha,\beta}e_{\alpha+\beta}$
for $c_{\alpha,\beta}\in F$, and $R(I) = 0$. By Lemma~5, $R$ gives
an RB-operator of weight~0 on $L$.

{\bf Statement 4}.\it\
a$)$ Let $D\in\mathrm{Der}_F(S,I,A)$, $f(x,y) = D(xy)-D(x)y-xD(y)$.
Assume that there exists
$\theta\colon S\to I$ such that
$-f(x,y) = \theta(xy)-\theta(x)y-x\theta(y)$
for all $x,y\in S$. Then $D+\theta\in\mathrm{Der}_F(S,I,A)$.

b$)$ Let $A$ be an  algebra over a field $F$, and let
$A = S\oplus I$ for some subalgebra $S$ of $A$ and an ideal $I$ of $A$.
Then $\mathrm{Der}_F(S,I,A) = \mathrm{Der}_F(S) + \mathrm{End}_F(S,I);$
i.e., every derivation $D$ from $S$ in $A$ modulo $I$ is a sum of
a~derivation $D_1\in\mathrm{Der}(S)$ and a linear map
$\theta\colon S\to I$, and conversely.\rm

{\it  Proof} of $a$) is straightforward.

$b$) Take $D\in \mathrm{Der}(S,I,A)$. Put $D_1=\pi\circ D$,
$D_1:S\mapsto S$, where $\pi$ is the projection on $S$; i.~e.,
$D_1(s)=\pi(D(s))\in S$ for all $s\in S$. Then $\pi\in
\mathrm{Hom}_F(A,S)$. Now, it suffices to put $\theta=D-D_1$.

The converse assertion is immediate.\qedd

\section{RB-Operators of Nonzero Weight}
\subsection{The Simple Jordan Algebra of Bilinear Form}

Let $J_{n+1}(f) = F\cdot1\oplus V$ be the direct sum of $F$ and a linear
$n$-dimensional space $V$, $n\geq2$,
and let $f$ be a nondegenerate symmetric bilinear form on $V$.
With respect to the product
\begin{equation}\label{FormProduct}
(\alpha\cdot1 + a)(\beta\cdot1 + b) =
(\alpha\beta+f(a,b))\cdot1 + (\alpha b +\beta a),\quad
\alpha,\beta\in F,\ a,b\in V,
\end{equation}
the space $J_{n+1}(f)$ is a simple Jordan algebra \cite{Nearly}.

The algebra $J_{n+1}(f)$ is quadratic, since for every
$x = \alpha\cdot1+a\in J_{n+1}(f)$, $\alpha\in F$, $a\in V$, we have
$x^2 - 2\alpha x + (\alpha^2-f(a,a))\cdot1 = 0$.
Hence, $t(x) = 2\alpha$, $n(x) = \alpha^2-f(a,a)$.

Choose a basis $e_1$, $e_2$, \ldots, $e_n$ for $V$ such that
the matrix of the form $f$ in this basis is diagonal with
some elements $d_1,d_2,\ldots,d_n\in F$ on the main diagonal.
Since $f$ is nondegenerate, $d_i\neq0$ for each $i$.

Given an RB-operator $R$ of weight $\Delta$ on $J_{n+1}(f)$,
assume that  $R$ is defined by a matrix $(r_{ij})_{i,j=0}^n$
in the basis $1,e_1,e_2,\ldots,e_n$.

The identity \eqref{RB} is equivalent to the system of equations,
which is quadratic with respect to  $r_{ij}$.
Due to the symmetricity of $f$, it suffices to consider the equations
arising from the equalities by the products $x_0 y_0$, $x_s y_s$
(let us denote them as $\underline{00}$ and $\underline{ss}$
for $s>0$, respectively)
and $x_0 y_k + x_k y_0$, $x_k y_l + x_l y_k$, $k\neq l$
(notation: $\underline{0k}$ for $k>0$ and $\underline{kl}$ for $k\neq l$, $k,l>0$).
There are eight series of equations: (the bold number denotes
the projection of \eqref{RB} either on $1$ or on $e_i$)
\begin{eqnarray*}\allowdisplaybreaks
{\bf 0}, & \underline{00}: & d_1 r_{10}^2+\ldots+d_n r_{n0}^2 =
r_{00}^2+\Delta r_{00}+2(r_{01}r_{10}+\ldots+r_{0n}r_{n0}), \\
& \underline{ss}: & d_1 r_{1s}^2+\ldots+d_n r_{ns}^2 =
r_{0s}^2+d_s r_{00}(2r_{ss}+\Delta), \\
& \underline{0k}: & d_1 r_{10}r_{1k}+\ldots+d_n r_{n0}r_{nk}
 = \Delta r_{0k}+d_k r_{00}r_{k0} + r_{00}r_{0k}+\ldots+r_{0n}r_{nk}, \\
& \underline{kl}: & d_1 r_{1k}r_{1l}+\ldots+d_n r_{nk}r_{nl} =
r_{0k}r_{0l}+r_{00}(d_k r_{kl}+d_l r_{lk}), \\
{\bf i>0}, & \underline{00}: & 2(r_{i1}r_{10}+
\ldots+r_{in}r_{n0})+\Delta r_{i0}=0, \\
& \underline{ss}: & r_{i0}(2r_{ss}+\Delta)=0, \\
& \underline{1k}: & d_k r_{i0}r_{k0}+r_{i1}r_{1k}+
\ldots+r_{in}r_{nk}+\Delta r_{ik}=0, \\
& \underline{kl}: & r_{i0}(d_k r_{kl} + d_l r_{lk}) = 0.
\end{eqnarray*}

Assume that $R$ is  an RB-operator on $J_{n+1}(f)$ such that $R(1)\not\in F$
and $F$ is algebraically closed. So, we have
\begin{gather}
r_{ss} = -\Delta/2,\ s>0, \label{RB-systemA} \\
d_k r_{kl} + d_l r_{lk} = 0,\ k,l>0,k\neq l. \label{RB-systemB}
\end{gather}
Then the system of quadratic equations written above
is equivalent to the following
\begin{gather}
d_1 r_{10}^2+\ldots+d_n r_{n0}^2 =
r_{00}^2+\Delta r_{00}+2(r_{01}r_{10}+\ldots+r_{0n}r_{n0}), \label{RB-system1} \\
d_1 r_{1s}^2+\ldots+d_n r_{ns}^2 = r_{0s}^2,\ s>0, \label{RB-system2} \\
d_1 r_{10}r_{1k}+\ldots+d_n r_{n0}r_{nk}
 = \Delta r_{0k}+d_k r_{00}r_{k0} + r_{00}r_{0k}+
 \ldots+r_{0n}r_{nk},\ k>0, \label{RB-system3} \\
d_1 r_{1k}r_{1l}+\ldots+d_n r_{nk}r_{nl}
 = r_{0k}r_{0l},\ k,l>0,\ k\neq l, \label{RB-system4} \\
2(r_{i1}r_{10}+\ldots+r_{in}r_{n0})+\Delta r_{i0}=0,\ i>0, \label{RB-system5} \\
d_k r_{i0}r_{k0}+r_{i1}r_{1k}+\ldots+r_{in}r_{nk}+\Delta r_{ik}=0,\ i,k>0.
\label{RB-system6}
\end{gather}

By \eqref{RB-systemA} and \eqref{RB-system6}
for $i = k = s > 0$ by \eqref{RB-system2} we have
$$
r_{0s}^2
 = \sum\limits_{i=1}^n d_i r_{is}^2
 = -\sum\limits_{i=1}^n d_s r_{is}r_{si} + d_s\frac{\Delta^2}{2}
 = d_s(d_s r_{s0}^2 + \Delta r_{ss}) + d_s\frac{\Delta^2}{2}
 = d_s^2 r_{s0}^2.
$$
So,
\begin{equation} \label{RB-systemC}
r_{0s} = z_s d_s r_{s0}
\end{equation}
with $z_s\in\{-1,+1\}$. Therefore, \eqref{RB-system6} could be derived
from \eqref{RB-system2} with the help of \eqref{RB-systemA} and \eqref{RB-systemC}.

By \eqref{RB-systemA}, \eqref{RB-systemB}, and \eqref{RB-system4} we have
\begin{multline*}
r_{0k}r_{0l} = \sum\limits_{i=1}^n d_i r_{ik}r_{il}
 = -d_k\sum\limits_{i=1}^n r_{ki}r_{il} + 2d_k r_{kk}r_{kl} \\
 = d_k(\Delta r_{kl} + d_l r_{0k}r_{l0}) + 2d_k r_{kk}r_{kl}
 = d_k d_l r_{k0}r_{l0},
\end{multline*}
whence $z_s = z\in\{-1,+1\}$ for all $s>0$ by \eqref{RB-systemC}.

Applying \eqref{RB-systemA}, \eqref{RB-systemB}, \eqref{RB-system5},
we get from \eqref{RB-system3}
\begin{multline*}
\Delta r_{0k}+d_k r_{00}r_{k0} + r_{00}r_{0k}
 = \sum\limits_{i=1}^n d_i r_{i0}r_{ik} - \sum\limits_{i=1}^n r_{0i}r_{ik} \\
 = (1-z)\sum\limits_{i=1}^n d_i r_{i0}r_{ik}
 = -d_k(1-z)\left( \sum\limits_{i=1}^n r_{i0}r_{ki} \right)
   +(1-z)d_k(r_{k0}r_{kk}+r_{k0}r_{kk}) \\
 = (1/2)d_k(1-z)\Delta r_{k0} -d_k\Delta(1-z)r_{k0} = (1/2)d_k(z-1)\Delta r_{k0}.
\end{multline*}
Thus, $(1+z)r_{k0}(2r_{00}+\Delta) = 0$.
Since $R(1)\neq F$, $(1+z)(2r_{00}+\Delta) = 0$.

Summarizing, we have the following system on
$\bar{r}_{ij} = \frac{\sqrt{d_i}}{\sqrt{d_j}}r_{ij}$ satisfying
$\bar{r}_{kk} = -\Delta/2$, $\bar{r}_{kl} = -\bar{r}_{lk}$,
$\bar{r}_{0k} = z\bar{r}_{k0}$
for $k,l>0$, $k\neq l$, $z\in\{-1,+1\}$:
\begin{gather}
(1+z)(2\bar{r}_{00}+\Delta) = 0, \label{RB-system0new} \\
(1-2z)\sum\limits_{p=1}^n \bar{r}_{p0}^2 = \bar{r}_{00}(\bar{r}_{00}+\Delta),\quad
\sum\limits_{p=1}^n \bar{r}_{pi}\bar{r}_{p0}
 = -\frac{\Delta z}{2}\bar{r}_{0i},\ i>0,\label{RB-system1new} \\
\sum\limits_{p=1}^n \bar{r}_{pk}\bar{r}_{pl}
 = \bar{r}_{0k}\bar{r}_{0l},\ k,l>0, k\neq l,\quad
\sum\limits_{p=1}^n \bar{r}_{pk}^2 = \bar{r}_{0k}^2,\ k>0. \label{RB-system2new}
\end{gather}

Consider the first case: (I) $z = 1$, $\bar{r}_{00} = -\Delta/2$.
Then the system \eqref{RB-system1new}--\eqref{RB-system2new} is of the form
\begin{gather}
\sum\limits_{p=1}^n \bar{r}_{p0}^2 = \frac{\Delta^2}{4},\quad
\sum\limits_{p=0}^n \bar{r}_{pi}\bar{r}_{p0}
 = -\Delta\bar{r}_{0i},\ i>0,\label{RB-system1newI} \\
\sum\limits_{p=1}^n \bar{r}_{pk}\bar{r}_{pl}
 = \bar{r}_{0k}\bar{r}_{0l},\ k,l>0, k\neq l, \quad
\sum\limits_{p=1}^n \bar{r}_{pk}^2 = \bar{r}_{0k}^2,\ k>0. \label{RB-system2newI}
\end{gather}

The second case is the following: (II) $z = -1$
(in what follows we assume that $\char F\neq3$),
\begin{gather}
\sum\limits_{p=1}^n \bar{r}_{p0}^2
 = \frac{\bar{r}_{00}(\bar{r}_{00}+\Delta)}{3}, \quad
\sum\limits_{p=1}^n \bar{r}_{pi}\bar{r}_{p0}
 = \frac{\Delta}{2}\bar{r}_{0i},\ i>0, \label{RB-system1newII} \\
\sum\limits_{p=1}^n \bar{r}_{pk}\bar{r}_{pl}
 = \bar{r}_{0k}\bar{r}_{0l},\ k,l>0, k\neq l, \quad
\sum\limits_{p=1}^n \bar{r}_{pk}^2 = \bar{r}_{0k}^2,\ k>0. \label{RB-system2newII}
\end{gather}

Applying  \eqref{RB-system1newII}--\eqref{RB-system2newII}, we obtain
\begin{multline*}\allowdisplaybreaks
\sum\limits_{i=1}^n \bar{r}_{i0}^2
 = \frac{4}{\Delta^2}\sum\limits_{i=1}^n\left(\sum\limits_{p=1}^n
 \bar{r}_{pi}\bar{r}_{p0}\right)^2 \\
 = \frac{4}{\Delta^2}\sum\limits_{p=1}^n
 \bar{r}_{p0}^2\left(\sum\limits_{i=1}^n\bar{r}_{pi}^2\right)
 + \frac{8}{\Delta^2}\sum\limits_{p,q=1,\,p\neq q}^n \bar{r}_{p0}\bar{r}_{q0}
   \left(\sum\limits_{i=1}^n\bar{r}_{pi}\bar{r}_{qi}\right) \\
 = \frac{4}{\Delta^2}\sum\limits_{p=1}^n \bar{r}_{p0}^4
 + \frac{8}{\Delta^2}\sum\limits_{p,q=1,\,p\neq q}^n \bar{r}_{p0}^2\bar{r}_{q0}^2
 = \frac{4}{\Delta^2}\left(\sum\limits_{p=1}^n \bar{r}_{p0}^2 \right)^2,
\end{multline*}
whence $A = \sum\limits_{i=1}^n \bar{r}_{i0}^2$ is equal to 0 or $\Delta^2/4$.

Suppose that $A = 0$. By \eqref{RB-system1newII}, up to action of $\phi$
we may assume that $\bar{r}_{00} = 0$.
So, $R(1)R(1) = 0$. By \eqref{RB}, we have
\begin{equation}\label{Nilpotency}
\begin{gathered}
0 = R(1)R(1) = 2R^2(1) + \Delta R(1), \\
0 = R(1)R(1)R(1) = 2R^2(1)R(1) = 2R^3(1) + 2\Delta R^2(1)
 = \Delta R^2(1),
\end{gathered}
\end{equation}
whence $R^2(1) = 0 = R(1)$, a contradiction to the assumption $R(1)\not\in F$.

Thus, $A = \Delta^2/4$ and by \eqref{RB-system1newII} we arrive at the following subcases.

(II\,a) $z = -1$, $\bar{r}_{00} = \Delta/2$.
In this case, the system
\eqref{RB-system1newII}--\eqref{RB-system2newII} is of the form
\begin{gather}
\sum\limits_{p=1}^n \bar{r}_{p0}^2 = \frac{\Delta^2}{4},\quad
\sum\limits_{p=0}^n \bar{r}_{pi}\bar{r}_{p0}
 = \Delta\bar{r}_{0i},\ i>0, \label{RB-system1newIIa} \\
\sum\limits_{p=1}^n \bar{r}_{pk}\bar{r}_{pl}
 = \bar{r}_{0k}\bar{r}_{0l},\ k,l>0, k\neq l, \quad
\sum\limits_{p=1}^n \bar{r}_{pk}^2
 = \bar{r}_{0k}^2,\ k>0. \label{RB-system2newIIa}
\end{gather}

(II\,b) $z = -1$, $\bar{r}_{00} = -3\Delta/2$,
\begin{gather}
\sum\limits_{p=1}^n \bar{r}_{p0}^2 = \frac{\Delta^2}{4},\quad
\sum\limits_{p=0}^n \bar{r}_{pi}\bar{r}_{p0}
 = -\Delta\bar{r}_{0i},\ i>0, \label{RB-system1newIIb} \\
\sum\limits_{p=1}^n \bar{r}_{pk}\bar{r}_{pl}
 = \bar{r}_{0k}\bar{r}_{0l},\ k,l>0, k\neq l, \quad
\sum\limits_{p=1}^n \bar{r}_{pk}^2 = \bar{r}_{0k}^2,\ k>0. \label{RB-system2newIIb}
\end{gather}

Notice that the numbers $\bar{r}_{ij}$ satisfying (I) and (II\,a)
could be obtained from each other by multiplying the first row by $-1$.
Further, for both cases (II\,a) and (II\,b) we define
$$
s_{kl} = \begin{cases}
\Delta/2, & k = l, \\
\bar{r}_{kl}, & k,l>0,\,k\neq l, \\
i\bar{r}_{kl}, & k\ or\ l=0,\,k\neq l,
\end{cases}
$$
where $i$ is a root of $x^2 + 1 = 0$.

It is easy to prove that the systems
\eqref{RB-system1newIIa}--\eqref{RB-system2newIIa} and
\eqref{RB-system1newIIb}--\eqref{RB-system2newIIb}
in the terms of $s_{ij}$ have the same form
\begin{equation}\label{RB-systemFinal}
\sum\limits_{p=0}^n s_{pk}s_{pl} = 0,\quad 0\leq k,l\leq n.
\end{equation}
We can represent the matrix $S = \{s_{kl}\}$ as $S = \frac{\Delta}{2}E + M$
for the skew-symmetric matrix~$M$ and the identity matrix $E$.

The system \eqref{RB-systemFinal} is equivalent to the equality
$(\frac{\Delta}{2}E+M)^T(\frac{\Delta}{2}E+M) = 0$ or,
applying the skew-symmetricity of $M$, we have $M^2 = \frac{\Delta^2}{4}E$.

{\bf Theorem 2}.\it \
Let $J_{n+1}(f)$ be the simple Jordan algebra of bilinear form $f$.
If $n$ is even then all RB-operators on $J_{n+1}(f)$ of nonzero weight
are splitting.\rm

{\it  Proof}.
Let $R$ be a non-splitting RB-operator of weight $\Delta\neq0$,
which is defined by a matrix $(r_{ij})_{i,j=0}^n$ in a basis $1,e_1,\ldots,e_n$.
By Lemma~1\,b), $R(1)\not\in F$. Let $\bar{F}$ be an algebraical closure of $F$.

Assume that $\char F\neq3$.
Then as it was stated above we can construct a skew-symmetric matrix
$M\in M_{n+1}(\bar{F})$ such that $M^2 = \frac{\Delta^2}{4}E$.
Hence, the rank of $M$ is equal to $n+1$.
It is well-known that the rank of a skew-symmetric matrix over the field
of characteristic different from 2 is even \cite{Gant}. We arrive at a contradiction.

If $\char F = 3$ then in the case (II) we have the following system of equations:
\begin{gather}
\bar{r}_{00}(\bar{r}_{00}+\Delta) = 0, \quad
\sum\limits_{p=1}^n \bar{r}_{pi}\bar{r}_{p0}
 = \frac{\Delta}{2}\bar{r}_{0i},\ i>0, \label{RB-system1newII-3} \\
\sum\limits_{p=1}^n \bar{r}_{pk}\bar{r}_{pl}
 = \bar{r}_{0k}\bar{r}_{0l},\ k,l>0, k\neq l, \quad
\sum\limits_{p=1}^n \bar{r}_{pk}^2 = \bar{r}_{0k}^2,\ k>0.
\label{RB-system2newII-3}
\end{gather}
Up to action of $\phi$, we may assume that $\bar{r}_{00} = 0$. By
the same reasons as above, from
\eqref{RB-system1newII-3}--\eqref{RB-system2newII-3} we see that $A
= \sum\limits_{i=1}^n \bar{r}_{i0}^2$ is equal to 0 or $\Delta^2/4$.
As it was proved above, the case $A = 0$ is contradictory. For $A =
\Delta^2/4$, we define the matrix $Q = (q_{kl})\in M_{n+1}(\bar{F})$
with the entries
$$
q_{kl} = \begin{cases}
-\Delta/2, & k = l, \\
\bar{r}_{kl}, & k,l>0,\,k\neq l, \\
i\bar{r}_{kl}, & k\ or\ l=0,\,k\neq l.
\end{cases}
$$
Analogously, we obtain  $Q^T Q = 0$ and $Q = -\frac{\Delta}{2}E + M$
for a skew-symmetric matrix~$M$.
The final arguments are the same as in the case $\char F\neq3$.\qed

Actually, we have proved even more than Theorem 2 states:

{\bf Corollary 3}.\it\
Let $J_{n+1}(f)$ be the simple Jordan algebra of bilinear form $f$,
and let $R$ be an RB-operator on $J_{n+1}(f)$ of nonzero weight.
If $n$ is even then we have $R(1) = 0$ up to $\phi$.\rm

{\bf Remark 4}.
Notice that for the simple Jordan algebra $J_{n+1}(f)$
of bilinear form $f$ and odd $n$,
there is the correspondence between
the set $X_\Delta$ of all RB-operators of nonzero weight $\Delta$ on $J_{n+1}(f)$
with the property $R(1)\not\in F$ for all $R\in X$
and the set $Y_\Delta$ of all skew-symmetric matrices from $M_{n+1}(F)$ satisfying
$S^2 = \frac{\Delta^2}{4}E$ for $S\in Y_\Delta$.
It is interesting to compare with the weight zero case.
In \cite{Gub2017} it was proved that over an algebraically closed field $F$,
we have the correspondence between the set $X_0$ of RB-operators of weight zero on $J_{n+1}(f)$
satisfying $R(1)\not\in F$ and $R^2 = 0$ for $R\in Z$
and the set $Y_0$ of all skew-symmetric matrices from $M_{n+1}(F)$
whose squares are zero.

The following example says about the situation in even dimension
over an algebraically closed field.

{\bf Example 10}.
Let $J_{2n}(f)$ be the simple Jordan algebra of bilinear from $f$
over an algebraically closed field $F$.
The following operator $\frac{2}{\Delta}R$ defined
by nonzero matrix entries of $R$ as
\begin{gather*}
r_{00} = -3,\quad r_{01} = \sqrt{d_1},\quad r_{10} = -\frac{1}{\sqrt{d_1}},\\
r_{ii} = -1,\quad r_{i\,i+1} = \frac{d_{i+1}}{d_i}\sqrt{ -\frac{d_i}{d_{i+1}} },
\quad
r_{i+1\,i} = - \sqrt{ -\frac{d_i}{d_{i+1}} },\ i\geq1,
\end{gather*}
is a non-splitting RB-operator on $J_{2n}(f)$ of weight $\Delta$.
This RB-operator arises from the case (II\,b).

Example 10 may be generalized for the simple countable-dimensional Jordan
algebra of diagonalized bilinear form.

The next example shows that non-splitting RB-operators of nonzero weight
on the simple even-dimensional Jordan algebra of bilinear form
can not to be block diagonal (as in Example 10).

{\bf Example 11}.
Consider $J_4(f)$ over $\mathbb{Z}_5$ with the form $f$ having the identity matrix
in the basis $1,e_1,e_2,e_3$. Then the following operator on $J_4(f)$ (arisen from the case (II\,b))
\begin{gather*}
R(1) = 4 + 4e_1 + 3e_2 + 3e_3,\quad
R(e_1) = 1 + 3 e_1 + 4 e_2 + e_3,\\
R(e_2) = 2R(e_1),\quad
R(e_3) = 2 + 4 e_1 + 3 e_2 + 3e_3
\end{gather*}
is a non-splitting RB-operator of weight $-1$.

We can see that there are also splitting RB-operators
 using all RB-operators from the cases (I) or (II).

{\bf Example 12}.
Consider $J_4(f)$ over $\mathbb{Z}_{13}$ with the form $f$
 having the identity matrix
in the basis $1,e_1,e_2,e_3$.
 Then the following operator on $J_4(f)$ (arisen from the case (I))
$$
R(1) = R(e_1) = 7 + 7e_1 + 7e_2 + 9e_3,\quad
R(e_2) = 7 + 6 e_1 + 7 e_2 + 4 e_3,\
R(e_3) = 5R(e_2)
$$
is a splitting RB-operator of weight $-1$, although $R(1)\not\in F$.
Here we have $\ker R = \Span\{1-e_1,e_2-5e_3\}$ and
$\Imm R = \Span\{1+e_2,e_1+5e_3\}$.

{\bf Statement 5}.\it\
Let $A$ be the simple Jordan algebra of bilinear form,
and let $R$ be an RB-operator on $A$ of nonzero weight $\Delta$.
If $R(1) = 0$ then $\dim\ker R\geq2$.\rm

{\it  Proof}.
By Lemma~1\,b), $R$ is splitting. So, $1\in\ker R$ and $1\not\in\Imm R$.

Suppose that $\dim \ker R = 1$. From $0 = R(R(e_i)+\Delta e_i)$,
$i=1,\ldots,n-1$, we deduce that
$R(e_i) = r_i\cdot 1-\Delta e_i$ for all $i=1,\ldots,n-1$ and for some $r_i\in F$.
Since $R(e_1)R(e_2) = r_1r_2\cdot1 \in \Imm R$,
we obtain either $r_1 = 0$ or $r_2 = 0$. Taking $r_1=0$,
one has $R(e_1)R(e_1) = d_1\Delta^2\cdot1\in\Imm R$ with nonzero $d_1\in F$,
a contradiction. \qed

In \cite{Gub2017}, all RB-operators on $J_3(f)$ of weight zero were described.
We have very close result for nonzero weight.

{\bf Example 13}.
Let $J_3(f)$ be the simple 3-dimensional Jordan algebra of bilinear form $f = (d_1,d_2)$,
and let $R$ be a nontrivial RB-operator on $J_3(f)$ of nonzero weight~$\Delta$.
By Corollary 3, up to $\phi$ we have $R(1) = 0$.
By Statement 5, $\dim \Imm R = 1$.
Thus, $R(e_1) = k(\alpha_0\cdot1+\alpha_1 e_1+\alpha_2 e_2)$,
$R(e_2) = l(\alpha_0\cdot1+\alpha_1 x+\alpha_2 y)$
for some $k,l,\alpha_i\in F$, $k$ and $l$ are nonzero
simultaneously as well as $\alpha_i$.
We have $le_1-ke_2\in\ker R$, so $R$ is splitting with respect to the subalgebras
$A_1 = \langle 1,le_1-ke_2\rangle$
and $A_2 = \langle \alpha_0\cdot1+\alpha_1 e_1+\alpha_2 e_2\rangle$.
The image of $R$ is a subalgebra of $J_3(f)$, so
$\alpha_0^2 - d_1\alpha_1^2-d_2\alpha^2 = 0$.
By \eqref{RB}, $k\alpha_1+l\alpha_2+\Delta = 0$
(it corresponds to the fact that $J_3(f) = A_1\oplus A_2$).
Thus, we described all RB-operators on $A$
of nonzero weight up to $\phi$.

\subsection{(Anti)Commutator Algebras}

Given an algebra $A$ with a product $\cdot$,
define the  operations $\circ$ and $[\,,]$ on the vector space of
 $A$ by the rule
$$
a\circ b=a\cdot b+b\cdot a,\quad [a,b]=a\cdot b-b\cdot a.
$$
We denote the space $A$ with $\circ$ as $A^{(+)}$
and the space $A$ with $[\,,]$ as $A^{(-)}$.

{\bf Statement 6}.\it\
Given an RB-operator $R$  of weight $\Delta$ on an algebra $A$,
$R$ is an RB-operator on $A^{(+)}$ and $A^{(-)}$ of weight $\Delta$.\rm

{\it  Proof}
is immediate by \eqref{RB}.\qed

{\bf Corollary 4}.\it\
Given an algebra $A$, if all RB-operators on $A^{(+)}$ $($or $A^{(-)})$ of nonzero weight
are splitting, then all RB-operators on $A$ of nonzero weight are splitting.\rm

{\it  Proof}.
Let $R$ be an RB-operator of nonzero weight $\Delta$ on $A$.
By Statement 6, $R$ is an RB-operator of weight $\Delta\neq0$ on $A^{(+)}$ and $A^{(-)}$.
By hypothesis, $R(R+\Delta\id) = 0$ on  $A$.
Thus, $R$ is splitting on $A$ by Lemma~1\,a).\qedd

\subsection{The matrix algebra of order 2}

{\bf Example 14}.
Define a linear map $R$ on $M_n(F)$ as follows:
$R$ is zero on all strictly upper (lower) triangular
matrices; $R$ is equal to $-\id$ on all strictly lower (upper) triangular
matrices; $R$ is an RB-operator on the algebra of diagonal
matrices of weight 1 \cite{AnBai}.
Then $R$ is an RB-operator on $M_n(F)$ of weight 1.

For example, a linear map $R$ on $M_2(F)$ such that
$R(e_{11}) = R(e_{12}) = 0$, $R(e_{22}) = e_{11}$, and $R(e_{21}) = -e_{21}$
is an RB-operator on $M_2(F)$ of weight 1.

Due to \cite{AnBai}, the set of all RB-operators of Example 14
 is invariant under~$\phi$.

{\bf Lemma 6}.\it \
Let $A$ be a quadratic algebra with a unit 1, and let
$R$ be an RB-operator on $A$ of weight 1, which is non-splitting.
If $R(1) = \alpha\cdot1 + p$ with $t(p) = 0$
then one of three following cases occurs$:$

{\rm I)} $R(1) = -\frac{1}{2}+p$, $R(p) = \frac{1}{4}-\frac{p}{2};$

{\rm II)} $R(1) = \frac{1}{2}+p$, $R(p) = -\frac{1}{4}-\frac{p}{2};$

{\rm III)}~$R(1) = -\frac{3}{2}+p$, $R(p) = -\frac{1}{4}-\frac{p}{2}$.\rm

{\it  Proof}. By Lemma~1\,b), $p\not\in F$.
Let $R(p) = \Delta\cdot1+s$, where $t(s) = 0$. Then
\begin{multline}\label{MatCond1}
(\alpha^2-n(p))\cdot1+2\alpha p
 = R(1)R(1)
 = 2R(R(1))+R(1)\\
 = (2\alpha^2+\alpha)\cdot1 + (2\alpha+1)p + 2R(p).
\end{multline}
By \eqref{MatCond1}, we conclude
\begin{equation}\label{MatCond2}
\Delta = -\frac{1}{2}(n(p) + \alpha(\alpha+1)),\quad s = -\frac{p}{2}.
\end{equation}

Considering
\begin{multline}\label{MatCond3}\allowdisplaybreaks
\bigg(\alpha\Delta+\frac{1}{2}n(p)\bigg)\cdot1 + \bigg(\Delta-\frac{\alpha}2\bigg)p
 = (\alpha\cdot1+p)\bigg(\Delta\cdot1-\frac{p}2\bigg) \\
 = R(1)R(p)
 = R(R(1)p+R(p)+p) \\
 = R\bigg((\alpha\cdot1+p)p+\Delta\cdot1-\frac{p}2+p\bigg) \\
 = (\Delta-\det p)(\alpha\cdot1+p)
 + \bigg(\alpha+\frac{1}2\bigg)\bigg(\Delta\cdot1-\frac{p}2\bigg) \\
 = \bigg(2\alpha\Delta+\frac{\Delta}2-\alpha n(p)\bigg)\cdot1
  + \bigg(\Delta-\frac{\alpha}2-\frac{1}4- n(p)\bigg)p,
\end{multline}
we have $n(p) = -1/4$, and
$\big(\alpha+\frac{1}{2}\big)\big(\Delta+\frac{1}{4}\big) = 0$.
Solutions to the last equation give exactly the
required cases I, II, and III.

{\bf Theorem 3}. \it\
All RB-operators on $M_2(F)$ of nonzero weight either are splitting
or are defined by Example 14 up to conjugation by an automorphism of $M_2(F)$.
\rm

{\it  Proof}. Suppose that $R$ is an RB-operator on $M_2(F)$, which is non-splitting,
and $R(1) = \alpha\cdot1+p$, where $\alpha\in F$, $(0\neq)p\in \mathrm{sl}_2(F)$.
Apply Lemma~6. The case III is equivalent to the case II by $\phi(R) = -R-\id$.

Since $\det\big(\frac{1}{2}+p\big) = 0$ for $p\in\mathrm{sl}_2(F)$
and the square of $\big(\frac{1}{2}+p\big)$ is proportional
to itself in both cases I and II,
we can consider an RB-operator $P = R^{(\varphi)}$
with $\varphi\in\Aut(M_2(F))$ such that
$\varphi\big(\frac{1}{2}+p\big) = e_{11}$.
Hence, $P(e_{11}) = 0$. Let $P(e_{12}) = s$ and $P(e_{21}) = t$.

{\it  Case I}. $P(e_{22}) = -e_{22}$.
 We have
\begin{equation}\label{Mat2-1-I-1}
0 = P(e_{11})P(e_{12})
  = P(e_{11}s + e_{12}) = (1+s_{12})s.
\end{equation}
If $s = 0$ then
$$
0 = P(e_{21})P(e_{12})
  = P(te_{12}+e_{22})
  = t_{21}t-e_{22},
$$
a contradiction. Hence, $s_{12} = -1$.
Since $\Imm (P)$ is a subalgebra, $se_{22}\in \Imm (P)$
and $-e_{12}+s_{22}e_{22}\in\Imm(P)$.
Therefore, $e_{12}\in\Imm(P)$.

Consider
\begin{equation}\label{Mat2-1-I-2}
0 = P(e_{21})P(e_{11}) = (1+t_{21})t.
\end{equation}
If $t_{21} = -1$ then $e_{12}t = -e_{11}+t_{22}e_{12}\in\Imm(P)$
and $e_{11}\in\Imm(P)$. Further, $e_{21}\in\Imm(P)$
and $\Imm(P) = M_2(F)$, a contradiction.
Hence, $t = 0$. As $\dim\ker P = \dim\Imm(P) = 2$, we have
$s = -e_{12}+s_{22}e_{22}$. Comparing the expressions
$$
P(e_{12})P(e_{22})
 = (-e_{12}+s_{22}e_{22})(-e_{22}) = e_{12}-s_{22}e_{22},
$$
\vspace{-1cm}
\begin{multline*}
P( P(e_{12})e_{22} + e_{12}(P(e_{22}) + e_{22}) ) \\
 = P((-e_{12}+s_{22}e_{22})e_{22})
 = P(-e_{12}+s_{22}e_{22}) = e_{12} - 2s_{22}e_{22},
\end{multline*}
we have $s_{22} = 0$, and $P$ is splitting.

{\it  Case II}. $P(e_{22}) = e_{11}$.
Since $\mathrm{tr}(e_{12}) = \mathrm{tr}(e_{21}) = 0$,
$\det s = \det t = 0$ by Lemma~3\,a).
From
\begin{multline*}
P(e_{12})P(e_{22})
 = se_{11} = s_{11}e_{11}+s_{21}e_{21} \\
 = P(P(e_{12})e_{22}+e_{12}P(e_{22})+e_{12}e_{22})
 = P(se_{22}+e_{12})
 = s_{22}e_{11}+(1+s_{12})s
\end{multline*}
we see that $s = -e_{12}$ or $s = s_{11}e_{11}+s_{21}e_{21}$.
Analogously, considering $P(e_{22})P(e_{21})$,
we have either $t = -e_{21}$ or $t = t_{11}e_{11}+t_{12}e_{12}$.
Together with \eqref{Mat2-1-I-1} and \eqref{Mat2-1-I-2}
we have either $s = -e_{12}$ or $s = 0$, and either $t = -e_{21}$ or $t = 0$.
The case $s = -e_{12}$ and $t = -e_{21}$ leads to $\Imm (P) = M_2(F)$,
a contradiction. The case $s = t = 0$ leads to $e_{22} = e_{21}e_{12}\in\ker P$,
a contradiction. The cases $s = -e_{12}$, $t = 0$ and $s = 0$, $t = -e_{21}$
give the RB-operators from Example 14.\qedd

\subsection{The Grassmann algebra of plane}

Denote by $\mathrm{Gr}_2$ the Grassmann algebra of plane $\Span\{e_1,e_2\}$,
i.e., the elements $1,e_1,e_2,e_1\wedge e_2$ form
 a linear basis for $\mathrm{Gr}_2$.

The algebra $\mathrm{Gr}_2$ is quadratic, since for
$x = \alpha\cdot1+\beta e_1+\gamma e_2+\delta e_1\wedge e_2\in\mathrm{Gr}_2$
we have
$x^2 = \alpha^2\cdot1+2\alpha\beta e_1
 + 2\alpha\gamma e_2+2\alpha\delta e_1\wedge e_2
 = 2\alpha x-\alpha^2\cdot1$.
Hence, $t(x) = 2\alpha$, $n(x) = \alpha^2$.
Let $A_0 = \Span\{e_1,e_2,e_1\wedge e_2\}$.

{\bf Theorem 4}.\it\
All RB-operators of nonzero weight  on $\mathrm{Gr}_2$ are splitting.\rm

{\it  Proof}.
Suppose that $R$ is a non-splitting RB-operator of weight 1.
On the contrary, by Lemma~3, we have $n(R(x)) = 0$
for every $x\in A_0$. So, $t(R(x)) = 0$, $x\in A_0$.

Let $R(1) = \alpha\cdot1+p$, where $\alpha\in F$
 and $p$ is nonzero element in $A_0$.
By Lemma~6, we have $n(R(p))\neq0$ in all three cases, a contradiction.\qedd

\subsection{The simple Jordan superalgebra $\mathrm{K}_3$}

The simple Jordan superalgebra $\mathrm{K}_3$ is defined as follows:
$\mathrm{K}_3 = A_0\oplus A_1$, $A_0 = \Span\{e\}$ (the even part),
$A_1 = \Span\{x,y\}$ (the odd part),
$$
e^2 = e, \quad
ex = xe = \frac{x}{2},\quad
ey = ye = \frac{y}{2},\quad
xy = -yx = \frac{e}{2},\quad
x^2 = y^2 = 0.
$$

The superalgebra $\mathrm{K}_3$ is quadratic, because of $z^2 - t(z)z = 0$
for each $z\in \mathrm{K}_3$, and $t(\alpha e+\beta x+\gamma y) = \alpha$.

{\bf Theorem 5}.\it \
All RB-operators  of nonzero weight on $\mathrm{K}_3$ are splitting.\rm

{\it  Proof}.
Let $R$ be a non-splitting RB-operator on $\mathrm{K}_3$ of weight 1.
Applying \eqref{CircInQuad}, we have
\begin{multline}\label{RB:K3-1}
t(R(z))R(z) = R(z)R(z)
 = R(z\circ R(z)+z^2) \\
 = t(z)R(R(z))+t(R(z))R(z)+t(z)R(z),
\end{multline}
whence
$t(z)R(R(z)+z) = 0$ for all $z\in \mathrm{K}_3$.

Hence, $R(R(e)+e) = 0$ and
$R(R(e+s) +e+s)= R(R(s)+s) = 0$
for each $s\in A_1$. Combining the last two equalities
we obtain $R(R(z)+z) = 0$ for all $z\in \mathrm{K}_3$.
The statement follows by Lemma~1\,a).\qedd

\subsection{Derivations of Nonzero Weight}

Given an algebra $A$ and $\Delta\in F$,
a linear operator $d\colon A\rightarrow A$
is called a {\it derivation of weight} $\Delta$ \cite{DifOp}
 provided that the following equality holds for all $x,y \in A$:
\begin{equation}\label{Diff}
d(xy) = d(x)y + xd(y) + \Delta d(x)d(y).
\end{equation}

Let us call the zero operator and $-\Delta\id$
as {\it trivial derivations of weight} $\Delta$.

{\bf Statement 7}. \cite{preLie}\it\
Given an algebra $A$ and an invertible derivation $d$ on $A$ of weight~$\Delta$,
the operator $d^{-1}$ is an RB-operator on $A$ of weight $\Delta$.\rm

{\it  Proof}.
Let $x=d^{-1}(a)$ and $y=d^{-1}(b)$. Then
$d^{-1}$ acts on both sides of \eqref{Diff} by the rule:
$$\qquad\qquad\qquad\qquad\quad
d^{-1}(a)d^{-1}(b)=d^{-1}(ad^{-1}(b)+d^{-1}(a)b+\Delta ab).
\qquad\qquad\qquad\qquad\quad\ \square
$$

{\bf Corollary 5}.\it\
There are no nontrivial invertible derivations of nonzero weight on
 quadratic division algebras, the simple odd-dimensional Jordan algebras
of bilinear form, the matrix algebra $M_2(F)$, the Grassmann algebra
$\mathrm{Gr}_2$, and the Kaplansky superalgebra~$\mathrm{K}_3$.\rm

{\it  Proof} follows from Theorems 1--5.\qedd

\section{The RB-Operators of Weight Zero}

\subsection{The Matrix Algebra of Order 2}

{\bf Lemma 7}.\it\
Let $R$ be an RB-operator on $M_n(F)$ of weight zero,
and let $\char F = 0$. Then $\Imm R$ consists only of degenerate matrices,
and $\dim(\Imm R)\leq n^2-n$.\rm

{\it  Proof}.
If $\Imm R$ contains an invertible matrix, then $1\in\Imm R$
 by the Cayley-Hamilton theorem, a contradiction with Lemma~1\,c).
  Thus, by \cite{Meshulam} we have
$\dim (\Imm R)\leq n^2-n$.\qed

{\bf Theorem 6}. \cite{Mat2}\it\
All nonzero RB-operators of weight zero
on $M_2(F)$ over an algebraically closed field~$F$
up to conjugation by automorphisms of $M_2(F)$, transposition
and multiplication by a nonzero scalar are the following$:$

$(M1)$ $R(e_{21}) = e_{12}$, $R(e_{11}) = R(e_{12}) = R(e_{22}) = 0;$

$(M2)$ $R(e_{21}) = e_{11}$, $R(e_{11}) = R(e_{12}) = R(e_{22}) = 0;$

$(M3)$ $R(e_{21}) = e_{11}$, $R(e_{22}) = e_{12}$, $R(e_{11}) = R(e_{12}) = 0;$

$(M4)$ $R(e_{21}) =- e_{11}$, $R(e_{11}) = e_{12}$, $R(e_{12}) = R(e_{22}) = 0$.
\rm

{\it  Proof}.
Let $R$ be a nonzero RB-operator on $M_2(F)$ of weight zero.
By Lemma~1\,d) or Lemma~7, $\dim(\Imm R)$ is equal to 1 or 2.

Let $\dim(\Imm R) = 1$. If $\Imm R = \Span\{v\}$ then  $\det v = 0$ by Lemma~7.
A~Jordan form of $v$ is equal to $e_{11}$ or $e_{12}$. By Lemma~2,
$R^2 = 0$ in both cases.
Up to conjugation by an automorphism of $M_2(F)$, we may assume that
either $\Imm R = F\cdot e_{11}$ or $F\cdot e_{12}$.

Consider the case $\Imm R = F\cdot e_{12}$.
If $R(1) = \alpha e_{12}$ for $\alpha\in F^*$
then from $0 = R(1)R(x) = R(R(1)x+R(x)) = R(R(1)x)$
and $0 = R(xR(1))$ for $x = e_{21}$ we have
$R(e_{11}) = R(e_{22}) = 0$, a contradiction.
So, $R(1) = 0$. If $R(e_{11}) = k e_{12}\neq 0$ then
$0 = R(e_{11})R(e_{21}) = R( k e_{12}e_{21} + e_{11}R(e_{21})
 = k^2 e_{12}$, a contradiction. Thus,
$R(e_{11}) = R(e_{22}) = R(e_{12}) = 0$ and $R(e_{21}) = \alpha e_{12}$
for some $\alpha\in F^*$, and we arrive at (M1).

Let $\Imm R = F\cdot e_{11}$. If $R(1) = \alpha e_{11}$ for $\alpha\in F^*$
then considering $(1/\alpha)R(1)R(x) = e_{11}R(x) = R(e_{11}x)$ for $x = e_{22}$
we get $R(1) = R(e_{22}) = 0$.
If $R(e_{12}) = \alpha e_{11}$ and $R(e_{21}) = \beta e_{11}$
for $\alpha\beta\neq0$ then the equality
$\alpha\beta e_{11} = R(e_{12})R(e_{21})
= R(\alpha e_{11}e_{21} + \beta e_{12} e_{11}) = 0$
gives a contradiction. Hence, $R(e_{12}) = 0$, $R(e_{21}) = \alpha e_{11}$ or
$R(e_{21}) = 0$, $R(e_{12}) = \alpha e_{11}$ for some $\alpha\in F^*$, this is (M2).

Let $\dim(\Imm R) = 2$. If $\Imm R$ is nilpotent
 then up to conjugation by $\Aut(M_2(F))$,
we can consider $e_{12}\in\Imm R$ and nonzero
$x = \begin{pmatrix}
a & 0 \\
b & c
\end{pmatrix}\in \Imm R$. Since
$\mathrm{tr}(x) = \det (x) = 0$, we have $a = c = 0$
and $e_{21}\in\Imm R$. Thus, $e_{12}e_{21}\in\Imm R$,
 a contradiction with $\dim(\Imm R)=2$.

Therefore, $\Imm R$ contains an idempotent. Up to conjugation by $\Aut(M_2(F))$,
$e_{11}\in\Imm R$. Since $\Imm R$ is a subspace of $M_2(F)$ consisting only of
degenerate matrices of maximal possible dimension; therefore,
$\Imm R = \Span\{e_{11},e_{12}\}$ up to transposition by~\cite{Flanders}.
Assume that $R$ is not nilpotent, so $R(1) = \alpha e_{11} + \beta e_{12}\neq 0$.
If $\alpha = 0$ then $R(1)R(1) = 2R^2(1) = 0$, and we get $R(1) = 0$
by \eqref{Nilpotency}. For $\alpha\neq0$,  applying Lemma~2\,a)
 we arrive at a contradiction.
So, $R$ is nilpotent, and $\Imm R\cap \ker R \neq (0)$.

a) $\Imm R = \ker R$. Let
$R(x_0 = \alpha e_{21}+\beta e_{22}) = e_{11}$, and
$R(y_0 = \gamma e_{21}+\delta e_{22}) = e_{12}$. From
$R(x_0)R(x_0) = e_{11} = R(\alpha e_{21})$
we have $\beta = 0$. Thus, $R(e_{21}) = (1/\alpha)e_{11}$
and $\delta\neq0$.
Considering $R(x_0)R(y_0) = e_{12} = \alpha R(e_{22})$,
we conclude that $R(e_{22}) = (1/\alpha)e_{12}$ and $\gamma = 0$. This is (M3).

b) $\dim(\Imm R\cap \ker R) = 1$. Assume that
there exists $a\in \Imm R\cap \ker R$
such that $a^2 = \alpha a$ for $\alpha \in F^*$.
Up to conjugation by $\Aut(M_2(F))$,
$a = e_{11}$. As above, $\Imm R = \Span\{e_{11},e_{12}\}$ up to transposition.
Let a nonzero $x = \beta e_{21}+\gamma e_{22}+\delta e_{12}$ belongs to
 $\ker R$. Since
$\ker R$ is an $\Imm R$-module, we get
$e_{11}x = \delta e_{12}\in\ker R$, and
$xe_{12} = \beta e_{22}\in\ker R$.
So, $\delta = 0$ and $e_{22}\in\ker R$.
Hence, $R(1) = 0$, and by Lemma~1\,c) $R^2 = 0$ and $\Imm(R)\subset \ker(R)$,
a contradiction.

Therefore, $\Imm R\cap \ker R$ is nilpotent, and it is equal to $F\cdot e_{12}$.
Let a nonzero $x_0 = \alpha e_{21}+\beta e_{22}+\gamma e_{11}$ belongs to $\ker R$.
From $e_{11}x_0 = \gamma e_{11}\in\ker R$ we have $\gamma = 0$.
If $\alpha\neq 0$ then $x_0 e_{12} = \alpha e_{22}\in\ker R$.
Hence, $e_{12},e_{22}\in\ker R$, a contradiction.
Thus, $\alpha = 0$ and $e_{22}\in\ker R$.
Let $R(z_0 := \alpha e_{11}+\beta e_{21}) = e_{11}$, and
$R(t_0 := \gamma e_{11}+\delta e_{21}) = e_{12}$.
From $e_{11} = R(z_0)R(z_0) = R(2\alpha e_{11}+\beta e_{21})$
we obtain $\alpha = 0$. From
$e_{12} = R(z_0)R(t_0) = R(\gamma e_{11}+\beta e_{22}) = R(\gamma e_{11})$
we have $\delta = 0$. Finally,
$0 = R(e_{11})R(e_{21}) = R((1/\gamma)e_{11}+(1/\beta)e_{11})$,
whence $\gamma = -\beta$, and we arrive at (M4).\qed

{\bf Corollary 6}.\it\
The set of all RB-operators of weight zero on an $n$-dimen\-sional algebra $A$
up to conjugation by automorphisms of $A$ and multiplications
on nonzero scalars may be considered as a projective variety
$RB(A)$ in $\mathbb{P}^{n^2-1}$ defined by $n^3$ relations obtained
from~\eqref{RB}, which is written on a linear basis of $A$.
Thus, by Theorem~6, $RB(M_2(F))$ has four fixed points
under the action by conjugation by an $($anti$)$automorphism.
Indeed, $(M4)$ is the only one that doesn't satisfy $R^2 = 0$.
Further, $(M1)$ and $(M2)$ but not $(M3)$ satisfy the condition
that all minors of order 2  are zero in the image.
Finally, the image of $(M1)$ in $M_2(\mathbb{C})$ is abelian,
but one of $(M2)$ is not abelian. Thus, the corresponding linear and quadratic
relations distinguish $(M1)$ and $(M2)$.\qedd

\subsection{The Grassmann Algebra of Plane}

{\bf Proposition 1}.\it\
Up to conjugation by an automorphism of $\mathrm{Gr}_2$
an arbitrary RB-operator $R$ of weight zero on $\mathrm{Gr}_{2}$
with a linear basis $1,e_1,e_2,e_1\wedge e_2$ is the following one$:$
$R(1),R(e_1)\in\Span\{e_2,e_1\wedge e_2\}$,
$R(e_2) = R(e_1\wedge e_2) = 0$.\rm

{\it  Proof}.
a) Take $x = \alpha\cdot1+x'\in R(\mathrm{Gr}_2)$,
where $x'\in \Span\{e_1,e_2,e_1\wedge e_2\}$. Then
$(x-\alpha\cdot1)^2 = 0$. Since $R(\mathrm{Gr}_2)$
is a subalgebra of $\mathrm{Gr}_2$, $\alpha^n\cdot1\in R(\mathrm{Gr}_2)$.
By Lemma~1\,c), $\alpha = 0$.

Given $x = \alpha\cdot1+x'$ with $x'\in \Span\{e_2,e_1\wedge e_2\}$, we have
$$
0 = R(x)R(x) = R(R(x)x + xR(x))
  = R(R(x)x'+x'R(x)) + 2\alpha R(R(x))
  = 2\alpha R(R(x)).
$$
At first $R(R(1)) = 0$; at second $R(R(x)) = 0$ for all $x\in \mathrm{Gr}_{2}$.
Hence, $\Imm R\subset \ker R$, and $\dim(\Imm R)\leq2$.

Assume that there exist $x$ and $y$  such that
$x_1e_1 + x_2e_2$ and $y_1e_1 + y_2e_2$
are linearly independent,
$R(x) = x_1e_1 + x_2e_2 + x_{12}e_1\wedge e_2$, and
$R(y) = y_1e_1 + y_2e_2 + y_{12}e_1\wedge e_2$. By~\eqref{RB},
$R(x)R(y) = (x_1y_2-x_2y_1)e_1\wedge e_2\in R(\mathrm{Gr}_2)$.
From here $R(\mathrm{Gr}_2) = \Span\{e_1,e_2,e_1\wedge e_2\}$, which
contradicts to the fact that $\dim(\Imm R)\leq2$.
Therefore, $R(\mathrm{Gr}_2)$ is an abelian algebra.

Show that $e_1\wedge e_2\in\ker(R)$. Otherwise, there is $x\in\mathrm{Gr}_2$
that $R(e_1\wedge e_2)x = e_1\wedge e_2$.
Hence,
$$
0 = R(e_1\wedge e_2)R(x) = R(R(e_1\wedge e_2)x + e_1\wedge e_2R(x))
  = R(e_1\wedge e_2),
$$
a contradiction.

Let $R(e_1) = k\alpha e_1+k\beta e_2+\gamma e_1\wedge e_2$, and
$R(e_2) = l\alpha e_1+l\beta e_2+\delta e_1\wedge e_2$. Then
$$
0 = R(e_1)R(1)
  = R(R(e_1)+e_1R(1))
  = R(R(e_1))
  = \alpha R(e_1) + \beta R(e_2),
$$
whence $R(e_1)$ and $R(e_2)$ are linearly dependent.
Up to conjugation by an automorphism of $\mathrm{Gr}_2$,
we may assume that $R(e_2) = 0$ and $R(e_1)\in\Span\{e_2,e_1\wedge e_2\}$.
It is immediate that a linear map $R$ such that
$R(1),R(e_1)\in\Span\{e_2,e_1\wedge e_2\}\subseteq \ker R$
is an RB-operator on~$\mathrm{Gr}_2$.\qedd

\subsection{The Simple Jordan Superalgebra $\mathrm{K}_3$}

{\bf Proposition 2}.\it\
An arbitrary  RB-operator $R$ of weight zero on $\mathrm{K}_3$
up to conjugation by $\Aut(\mathrm{K}_3)$ is the following one$:$
$$
R(e) = R(x) = 0, \quad
R(y) = ae + bx,\ a,b\in F.
$$\rm

{\it  Proof}.
Let $R$ be a nonzero RB-operator on $\mathrm{K}_3$ of weight zero.
By analogy with the proof of Theorem 5 and \eqref{RB:K3-1},
we have $R(R(z)) = 0$ for all $z\in \mathrm{K}_3$.
So, $\Imm R\subset \ker R$, and $\dim\Imm R = 1$, $\dim\ker R = 2$.

Let $R(e) = \alpha e+\beta x+\gamma y$. Suppose that $\alpha\neq0$.
For all $z\in \ker R$,
\begin{equation}\label{K3-0}
0 = R(e)R(z) = R(R(e)z)
\end{equation}
 hold, whence $R(e)z\in\ker R$. Since
 $R(e)\in\ker R$, there exists $z = \Delta x+\mu y\in\ker R$
 such that $z$ and $\beta x+\gamma y$ are linearly independent.
From $R(e)z = \alpha z+(\beta\mu-\gamma\Delta)e\in \ker R$,
we have $e\in \ker R$.

By analogy with \eqref{K3-0}, for $\alpha = 0$,
considering $z = \chi e+\Delta x+\mu y\in\ker R$,
we obtain $R(e) = 0$. From here and $\Imm R\subset \ker R$,
the assertion follows easily.\qedd

\subsection{Connection with the Yang---Baxter Equation}

Let $A$ be an associative algebra, $r = \sum a_i\otimes b_i\in
A\otimes A$. The tensor $r$ is called a solution of the associative
Yang---Baxter equation (AYBE, \cite{Zhelyabin,Aguiar01}) if
\begin{equation}\label{AYBE}
r_{13}r_{12}-r_{12}r_{23}+r_{23}r_{13} = 0,
\end{equation}
where
$$
r_{12} = \sum a_i\otimes b_i\otimes 1,\quad
r_{13} = \sum a_i\otimes 1\otimes b_i,\quad
r_{23} = \sum 1\otimes a_i\otimes b_i
$$
are elements from $A^{\otimes 3}$.

{\bf Example 15} \cite{Aguiar00}.
Let $r = \sum a_i\otimes b_i$ be a solution of AYBE on an associative algebra~$A$.
A linear map $R\colon A\to A$ defined as
$R(x) = \sum a_i x b_i$ is an RB-operator of weight zero on~$A$.

The image of an RB-operator of weight zero on an algebra $A$
is a subalgebra of $A$. The following example shows that the
kernel of an RB-operator of weight zero on $A$ is not a subalgebra of $A$
in general (even in the associative case).

{\bf Example 16}.
Consider the following solution to \eqref{AYBE} on $A = M_4(F)$ with an arbitrary field $F$:
$$
r = e_{11}\otimes e_{12}-e_{12}\otimes e_{11}
  + e_{33}\otimes e_{34}-e_{34}\otimes e_{33}.
$$
By Example 15,
$$
R(x) = e_{11}x e_{12}-e_{12}x e_{11} + e_{33}x e_{34}-e_{34}x e_{33}
$$
is an RB-operator on $A$, and its kernel consists of
the matrices $(a_{ij})$ in $A$ such that $a_{11}=a_{21}=a_{33}=a_{43}=0$.
It is easy to see that $\ker R$ is not a subalgebra of $A$.

Let $A$ be an algebra. Assume that $r = \sum a_i\otimes b_i\in
A\otimes A$ satisfies the {\it nonassociative Yang---Baxter equation}
over $A$ ($r_{12}=\sum a_i\otimes b_i\otimes 1, r_{13}=\sum
a_i\otimes 1 \otimes b_i$ and so on):
\begin{equation}\label{q6}
r_{12}r_{13}-r_{23}r_{12}-r_{13}r_{32}=0.
\end{equation}
(Note that the classical Yang---Baxter equation \cite{BelaDrin82} on
Lie algebras is written as
$[r_{12},r_{13}]+[r_{12},r_{23}]+[r_{13},r_{23}] = 0$.) We have
$$
\sum a_ia_j\otimes b_i\otimes b_j -a_i\otimes a_jb_i
\otimes b_j -a_i\otimes b_j\otimes b_ia_j=0,
$$
and if $\varphi\colon A\to A^*$ is one-to-one then
\begin{equation}\label{q7}
\sum a_ia_jb_i^*(x)b_j^*(y) -a_i(a_jb_i)^*(x)b_j^*(y)
 - a_ib_j^*(x)(b_ia_j)^*(y)=0.
\end{equation}

On the other hand, let $R(x)=\sum a_i b_i^*(x)$
be an RB-operator on $A$ of weight zero. Then \eqref{RB} gives
\begin{equation}\label{q8}
\sum a_ia_jb_i^*(x)b_j^*(y) -a_ib_j^*(y)b_{i}^{*}(xa_j)
-a_ib_j^*(x)b_{i}^{*}(a_jy)=0.
\end{equation}

Let $A$ be a finite-dimensional algebra such that the map
$\varphi\colon A\mapsto A^*$, $\varphi(a)=a^*$, is an
$A$-bimodule isomorphism, where the action on $A^*$
is defined by the rule
$$
a\cdot b^*(x)=b^*(xa),\quad b^*\cdot a(x)=b^*(ax).
$$
Let $R$ be an RB-operator on $A$ of weight zero.
If $a_1,\ldots,a_n$ is a basis for $A$ then
$R(a_i)=\sum\limits_{j=1}^n \alpha_{ij}a_j$. Define a linear functional
$b^*_j$ in $A^*$ by the rule
$b^*_j(a_i)=\alpha_{ij}$. Then
$R(a_i)=\sum\limits_{j=1}^n b^*_j(a_i)a_j$, whence
$R(x)=\sum\limits_{j=1}^n b^*_j(x)a_j$ for every $x\in A$.

If $\varphi\colon A\to A^*$, $\varphi(a)=a^*$, is an $A$-bimodule
isomorphism, where the action on $A^*$ is defined as above, then
it is easy to see that \eqref{q8} is equivalent to \eqref{q7}. Thus, we have

{\bf Theorem 7}.\it\
Let $A$ be a finite-dimensional algebra such that the map
$\varphi\colon A\to A^*$, $\varphi(a)=a^*$, is an $A$-bimodule
isomorphism, where the action on $A^*$ is defined by the rule
$$
a\cdot b^*(x)=b^*(xa),\quad b^*\cdot a(x)=b^*(ax).
$$
Then a solution $r=\sum a_i\otimes b_i$ of the nonassociative
Yang---Baxter equation \eqref{q6} over $A$ defines an RB-operator on
$A$ of weight zero by the rule $R(x)=\sum a_i b_i^*(x)$. Conversely,
if $R$ is an RB-operator on $A$ then $R(x)=\sum\limits_{j=1}^n
b^*_j(x)a_j$, and $r=\sum a_i\otimes b_i$ is a solution of the
nonassociative Yang---Baxter equation \eqref{q6} over $A$.\rm\qed

Let $A$ be a simple finite-dimensional algebra in some variety ${\cal M}$.
Assume that $A^*$ is an $A$-${\cal M}$-bimodule (with the action as above).
Then $A^*$ is an irreducible $A$-bimodule. Indeed, if $V$ is a submodule of $A^*$
and $V\neq A^*$ then there is $x\in A$ such that $f(x) = 0$ for all $f\in V$.
Then $a\cdot f(x)=0$ and $f\cdot a(x)=0$ for every $a\in A$, whence $f(I)=0$,
where $I = ideal\langle x\rangle = A$, i.e., $f = 0$.

{\bf Corollary 7}.\it\ Let $A$ be a simple finite-dimensional
self-adjoined algebra, i.e., $A\cong A^*$ as $A$-bimodules. Then the
solutions of the nonassociative Yang---Baxter equation $(\ref{q6})$
over $A$ are in one-to-one correspondence with the RB-operators on
$A$ of weight zero.\rm\qed

If $C$ is Cayley---Dickson algebra then $C$ is an alternative
$D$-bialgebra \cite{Goncharov}, and $C$ is a  self-adjoined algebra.
Therefore, we obtain the following

{\bf Corollary 8}.\it The solutions of the nonassociative
Yang---Baxter equation \eqref{q6} over $C$ are in one-to-one
correspondence with the RB-operators on $C$ of weight zero.\rm

Note that  all skew-symmetric ($r=\sum a_i\otimes b_i=-\sum
b_i\otimes a_i$) solutions of the Yang---Baxter equation over the
Cayley---Dickson matrix algebra were described in \cite{Goncharov}.

\section*{Acknowledgments}
Vsevolod Gubarev is supported by the Austrian Science Foun\-da\-tion FWF,
 grant P28079.


\begin{thebibliography}{99}

\bibitem{Aguiar00}
Aguiar, M.: Pre-Poisson algebras. Lett. Math. Phys.  {\bf 54},
263--277 (2000)

\bibitem{Aguiar01}
 Aguiar, M.:
On the Associative Analog of Lie Bialgebras. J. Algebra.  {\bf 244},
492--532 (2001)

\bibitem{AnBai}
An, H.,  Bai,  C.: From Rota---Baxter Algebras to Pre-Lie Algebras.
J. Phys. A.  N. 1, 015201. P. 19 (2008)

\bibitem{BBGN2011}
Bai, C.,  Bellier,  O., Guo, L., Ni, X.: Splitting of operations,
Manin products, and Rota---Baxter operators. Int. Math. Res. Notices.
 {\bf 3}, 485--524 (2013)

\bibitem{GB}
Baxter¡ G.: An analytic problem whose solution follows from a simple
algebraic identity. Pacific J. Math. {\bf 10}, 731--742 (1960)

\bibitem{BelaDrin82}
Belavin, A. A.,  Drinfel'd, V. G.: Solutions of the classical
Yang---Baxter equation for simple Lie algebras. Funct.~Anal.~Appl.
 {\bf 16}, 159--180 (1982)

\bibitem{Braga}
 de Bragan\c{c}a, S. L.:
Finite Dimensional Baxter Algebras. Studies in Applied Mathematics.
 {\bf ~54}(1), 75--89 (1975)

\bibitem{EF}
Ebrahimi-Fard, K.:  Loday-type algebras and the Rota---Baxter
relation.  Lett. Math. Phys. {\bf 61}(2),~139--147 (2002)

\bibitem{Flanders}
Flanders, H.: On spaces of linear transformations with bounded rank.
 J. London Math. Soc. {\bf 37},~10--16 (1962)

\bibitem{Gant}
Gantmacher, F. R.:  The theory of Matrices. V.~2. N.-Y.: Chelsea,
 (1959)

\bibitem{Goncharov}
Goncharov, M. E.: The classical Yang-Baxter equation on alternative
algebras: The alternative D-bialgebra structure on Cayley-Dickson
matrix algebras. Siberian Math. J. {\bf 48}(5),~809--823 (2007)

\bibitem{Goncharov1}
Goncharov, M. E.: On Rota-Baxter operators of non-zero weight arisen
from the solutions of the classical Yang---Baxter equation. Sib.
Electron. Math. Reports. {\bf 14}, 1533--1544 (2017)

\bibitem{Gub2017}
Gubarev, V.:   Rota---Baxter operators of weight zero on simple
Jordan algebra of Clifford type. Sib. Electron. Math. Reports. {\bf
 14}, 1524--1532 (2017)

\bibitem{Embedding}
Gubarev, V.,  Kolesnikov, P.: Embedding of dendriform algebras into
Rota---Baxter algebras. Cent. Eur. J. Math. {\bf 11}(2), 226--245
(2013)

\bibitem{GuoMonograph}
Guo, L.: An Introduction to Rota---Baxter Algebra. Surveys of Modern
Mathematics.  V.~4. Somerville, MA: International Press; Beijing:
Higher education press. (2012)

\bibitem{GEF}
Guo, L.,  Ebrahimi-Fard, K.: Rota---Baxter algebras and dendriform
algebras. J.~Pure Appl. Algebra. {\bf 212}, 320--339 (2008)

\bibitem{DifOp}
Guo, L.,  Keigher, W.: On differential Rota---Baxter algebras. J.
Pure Appl. Algebra.~{\bf 212}, 522--540 (2008)

\bibitem{Jian}
Jian, R.-Q.:  Quasi-idempotent Rota-Baxter operators arising from
quasi-idempotent elements. Lett. Math. Phys. {\bf 107}(2), 367--374
(2017)

\bibitem{Kolesnikov}
Kolesnikov, P. S.: Homogeneous averaging operators on simple finite
conformal Lie algebras. J. Math. Phys. {\bf 56}, 071702, (2015)

\bibitem{preLie}
Li, X. X.,  Hou, D. P., Bai, C.M.:  Rota-Baxter operators on pre-Lie
algebras. J.~Nonlinear Math. Phys. {\bf 14}(2), 269--289 (2007)

\bibitem{Meshulam}
Meshulam, R.: On the maximal rank in a subspace of matrices. Quart.
J. Math. Oxford. ~{\bf 36}(2), 225--229 (1985)

\bibitem{sl2}
Pan, Yu, Liu, Q., Bai, C., Guo,  L.: PostLie algebra structures on
the Lie algebra $\mathrm{sl}(2,\mathbb{C})$. Electron. J. Linear
Algebra. {\bf ~23},~180--197 (2012)

\bibitem{sl2-0}
Pei, J.,  Bai, C., Guo, L.: Rota-Baxter operators on
$\mathrm{sl}(2,\mathbb{C})$ and solutions of the classical
Yang-Baxter equation. J. Math. Phys. {\bf 55}, 021701, (2014)

\bibitem{Poj}
Pozhidaev, A. P.: 0-Dialgebras with bar-unity and nonassociative
Rota-Baxter algebras. Sib. Math. J. {\bf 50}(6), 1070--1080 (2009)

\bibitem{RG1}
Rota, G.-C.: Baxter algebras and combinatorial identities. I. Bull.
Amer. Math. Soc.~{\bf 75}, 325--329 (1969)

\bibitem{RG2}
Rota, G.-C.: Gian-Carlo Rota on combinatorics, introductory papers
and commentaries. Boston: Birkh\"{a}user, (1995)

\bibitem{Semenov83}
Semenov-Tyan-Shanskii, M. A.: What is a classical $r$-matrix? Funct.
Anal. Appl.~{\bf 17}, 259--272 (1983)

\bibitem{Mat2}
Tang, X., Zhang,  Y., Sun, Q.: Rota-Baxter operators on
4-dimensional complex simple associative algebras. Appl. Math. Comp.
 {\bf 229}, 173--186 (2014)

\bibitem{Zhelyabin}
Zhelyabin, V. N.: Jordan bialgebras of symmetric elements and Lie
bialgebras. Sib. Mat. J. {\bf 39}(2), 261--276 (1998)

\bibitem{Nearly}
Zhevlakov, K. A., Slin'ko, A. M., Shestakov, I. P., Shirshov, A. I.:
Rings that are nearly associative. N.-Y.: Academic Press, (1982)

\bigskip

Pilar Benito \\
Universidad de La Rioja\\
Calle Madre de Dios, 853, 26004\\
Logro\~{n}o, Spain \\
e-mail: {\sl pilar.benito@unirioja.es}\vspace*{1mm}

Vsevolod Gubarev \\
University of Vienna\\ 
Oskar-Morgenstern-Platz 1, 1090, Vienna, Austria \\
Sobolev Institute of mathematics\\
Acad. Koptyug ave., 4, 630090 Novosibirsk, Russia \\
e-mail: {\sl vsevolod.gubarev@univie.ac.at}\vspace*{1mm}

Alexander Pozhidaev \\
Sobolev Institute of mathematics \\
Acad. Koptyug ave., 4, 630090 Novosibirsk, Russia \\
Novosibirsk State University\\ 
Pirogova str., 2, 630090 Novosibirsk, Russia \\
e-mail: {\sl  app@math.nsc.ru}

\end{thebibliography}
\end{document}